\documentclass[12pt,a4paper]{amsart}
\usepackage{amssymb,mathrsfs}
\newcommand{\n}{\noindent}
\newcommand{\la}{\lambda}
\newcommand{\su}{\subset}
\newcommand{\m}{\medskip}
\newcommand{\ep}{\varepsilon}
\newcommand{\si}{\sigma}
\newcommand{\de}{\delta}

\newcommand{\al}{\alpha}

\newcommand{\p}{\partial}
\newcommand{\mc}{\mathcal}

\newcommand{\st}{($*$)}
\newcommand{\finproof}{\phantom{ii} \hfill $\blacksquare$ \bigskip}

\newtheorem{theorem}{Theorem}

\newtheorem{definition}[theorem]{Definition}
\newtheorem{fact}[theorem]{Fact}
\newtheorem{proposition}[theorem]{Proposition}
\newtheorem{corollary}[theorem]{Corollary}
\newtheorem{example}[theorem]{Example}

\begin{document}
\title{A note on fragmentability and weak-$G_\delta$ sets.}
\author{V. P. Fonf}\thanks{Department of Mathematics, Ben Gurion University of the Negev, Beer-Sheva, Israel ({\tt fonf@cs.bgu.ac.il})}
\author{R. J. Smith}\thanks{Institute of Mathematics of the AS CR, \v{Z}itn\'{a} 25, CZ - 115 67 Praha 1, Czech Republic ({\tt smith@math.cas.cz})}
\author{S. Troyanski}\thanks{Departamento de Matem\'aticas, Universidad de Murcia, Campus de Espinardo, 30100 Murcia, Spain ({\tt stroya@um.es})}
\thanks{V.\ P.\ Fonf was partially supported by the Spanish government, grant MEC SAB 2005-016. S.\ Troyanski was partially supported by the Center for Advanced
Studies at Ben-Gurion University of the Negev.}
\subjclass[2000]{}
\date{December 2008}
\begin{abstract}
In terms of fragmentability, we describe a new class of Banach
spaces which do not contain weak-$G_\delta$ open bounded subsets. In
particular, none of these spaces is isomorphic to a separable
polyhedral space.
\end{abstract}
\maketitle

\section{Introduction and Preliminaries}
All Banach spaces under consideration in this note are assumed to be
real and infinite-dimensional.

According to a well known theorem of Lindenstrauss and Phelps
\cite{LP}, if $X$ is a reflexive space then every closed convex
and bounded body in $X$ has uncountably many extreme points. The
first named author has obtained different generalisations of this
result. In particular, in \cite{F1}, it is proved that every
infinite-dimensional Banach space $X$ which is not
$c_0$-saturated does not admit a countable boundary. Moreover, if
$X$ is not $c_0$-saturated then \cite{F2}

\m

\n (a) {\rm $X$ does not contain an open, bounded weak-$G_\delta$
set.}

\m

\n In \cite{F3} it is shown that if a separable space $X$ does not
contain $c_0$ then

\m

\n (b) {\rm the polar $A^{\circ}$ of any closed convex and bounded
body $A\su X$ with $0\in {\rm int}A$ contains uncountably many
$w^*$-exposed points.}

\m

Recall that a set $B \subseteq B_{X^*}$ is said to be a boundary
for $X$ if, for every $x \in X$, there is $f \in B$ such that
$f(x) = ||x||$. Assume that $\{f_n\}_{n=1}^\infty$ is a countable
boundary for $X$ (the case: $X$ is polyhedral \cite{F4}), then
$$
{\rm int}\,B_X = \bigcap_{n=1}^\infty \{x \in X:\, f_n(x) < 1\}.
$$
Thus $\rm{int}\,B_X$ is an open, bounded weak-$G_\delta$ set.
Next let $G\su X$ be a bounded, convex, closed body and $0\in
{\rm int}G.$ A point $x \in \p G$ is said to be smooth if the
Minkowski functional $p$ of $G$ is G\^{a}teaux differentiable at
$x.$ A point $f$ of a subset $A\su X^*$ of the dual space $X^*$
is said to be a $w^*$-exposed point of $A$ if there is $x \in X$
such that $f(x)
> g(x)$ for every $g \in A$, $g \neq f$. Moreover, we say that
this $x$ $w^*$-exposes $f$. Let us recall also the following well
known fact.

\begin{fact}\label{ff}
A point $x \in G$ is smooth for $G$ if and only if $x$
$w^*$-exposes some point $f$ in the polar $G^{\circ}$ of $G$,
with $f(x) = 1$.
\end{fact}

\m

In this note we describe (by means of special fragmentable sets) a
new class $\mc{K}$ of Banach spaces $X$ which have both properties
(a) and (b) and which may be $c_0$-saturated.

\begin{definition}[Namioka {\cite{N}}]
A set $M$ in a Banach space $X$ is said to be {\em fragmentable} if,
for any subset $A$ of $M$ and any $\ep > 0$, there is a weak open
set $V$ which meets $A$ and ${\rm diam}(A \cap V) < \ep$.
\end{definition}

Additionally, a set $M$ is {\em dentable} if, for any $\ep > 0$,
there is an open half space $H$ which meets $M$ and ${\rm diam}(A
\cap H) < \ep$. Clearly, if every subset of a set $M$ is dentable
then $M$ is fragmentable. It is known that if $M$ is a weakly
compact subset of a Banach space, or a bounded subset of a dual
space of an Asplund space, then every subset of $M$ is dentable (see
e.g.\ \cite[p. 31, 60 and 91]{B}).

We define the class $\mc{K}$ as follows.

\begin{definition}\label{dk}
A Banach space $X$ belongs to $\mc{K}$ if $X$ contains a non-empty
fragmentable set $M \su {\rm int}B_X$ satisfying the following
condition

\begin{enumerate}
\item[\em{{\st}}] {\em for any $\ep > 0$, any weak open set $V$ and
any $x_0 \in V \cap M$, there is a finite sequence
$\{x_i\}_{i=1}^n \su V \cap M$ such that $||x_i -x_{i-1}|| <
\ep$, $i=1,\ldots,n,$ and $||x_n|| \geq 1 - \ep$.}
\end{enumerate}
\end{definition}

\m

Our main result is the following
\begin{theorem}\label{t1}
If $X\in\mc{K}$ then $X$ does not contain open bounded
$w-G_{\de}$ sets. Moreover, if $X\in\mc{K}$ is separable and $G$
is a convex bounded open set then the set of all smooth points of
${\rm cl}G$ cannot be covered by a countable union of weak closed
sets which does not meet $G$. In particular, if $0\in G$ then the
set $w^*$-exp $G^{\circ}$ is uncountable.
\end{theorem}

The following corollary complements the main result from
\cite{F3}.
\begin{corollary}
Assume that a separable Banach space $X$ contains a subspace with
the Radon-Nikod\'{y}m property (e.g. reflexive or $l_1$). If $F\su
X$ is a bounded closed convex body with $0\in {\rm int}F$ then the
set $w^*$-exp $F^{\circ}$ is uncountable.
\end{corollary}

Hitherto, there were no $c_0$-saturated Banach spaces known to
satisfy the conclusions of Theorem \ref{t1}. Despite the fact that
the proof of Theorem \ref{t1} uses some ideas from \cite{F3}, we
give examples of $c_0$-saturated Banach spaces which admit sets
satisfying the hypotheses of Theorem \ref{t1}.

\section{Proof of Theorem \ref{t1}}

The proof of the following fact is standard.

\begin{fact}\label{f1}Let $K$ be a weak compact subset of a
weak open subset $V$ of a Banach space $X.$ Then there is a
non-empty weak open neighbourhood $W$ of the origin such that $K + W
\su V$.
\end{fact}

The following rather technical proposition will be our main tool.

\begin{proposition}\label{p1}
Let $G$ be an open bounded subset of a Banach space $X$ and assume
$0 \in K \su G$, where $K$ is compact. Put $G_K = \{x \in X:\; x + K
\su G\}$. Assume that $M$ is a non-empty fragmentable subset of
${\rm cl}G_K$, such that for any weak open set $U$ with $U \cap M
\neq \emptyset$, for any weak closed subset $E$ with $E \cap G =
\emptyset$, and for any $\ep
> 0$, there is $y \in M \cap U$ such that $(y + K) \cap E =
\emptyset$ and $d(y, \p G_K) < \ep$. Then, for any $w$-$F_{\si}$ set
$F$ with $F \cap G = \emptyset$, there is $x \in X$ such that $x +
K\su {\rm cl}G$, $(x + K) \cap F = \emptyset$, and $(x + K) \cap \p
G \neq \emptyset$.
\end{proposition}

\n{\bf Proof.} Let $F = \bigcup_{n=1}^{\infty}F_n$, where $\{F_n\}$
is an increasing sequence of weak closed sets. Set $F_0 = \emptyset$
and let $\{\ep_n\}_{n=0}^\infty$ be a sequence of positive numbers
tending to $0$, where $\ep_0 > {\rm diam}(G_K)$. We construct a
sequence $\{x_n\}_{n=0}^\infty \su M$ and decreasing sequences of
$w$-open sets $\{U_n\}_{n=0}^\infty$ and $\{V_n\}_{n=0}^\infty$ with
the following properties
\begin{enumerate}
\item $x_n \in U_n$ and $x_n + K \su V_n$;
\item $w$-${\rm cl} V_n \cap F_n = \emptyset$;
\item ${\rm diam}(U_n \cap M) <\ep_n$;
\item $d(x_n ,\p G_K) < \ep_n$
\end{enumerate}
for all $n$. To begin, let $x_0 \in M$ be arbitrary and $U_0 = V_0 =
X$. Assume we have constructed $x_n$, $U_n$ and $V_n$. By Fact
\ref{f1}, we can take a weak open neighbourhood $W$ of $x_n$ such
that $W + K \su V_n$. Since $x_n \in U_n \cap W \cap M$ and $M$ is
fragmentable, there exists weak open $U_{n+1} \su W \cap U_n$ such
that $U_{n+1} \cap M$ is non-empty and ${\rm diam}(U_{n+1} \cap M) <
\ep_{n+1}$. From our hypothesis, there exists $x_{n+1} \in U_{n+1}
\cap M$ with the property that $(x_{n+1} + K) \cap F_{n+1} =
\emptyset$. Since $x_{n+1} + K \su U_{n+1} + K \su W + K \su V_n$
and $x_{n+1} + K \subseteq X\backslash F_{n+1}$, again by Fact
\ref{f1} we can pick a weak open neighbourhood $W^\prime$ of
$x_{n+1}$, satisfying $w$-${\rm cl}W^\prime + K \subseteq
V_n\backslash F_{n+1}$. Define $V_{n+1} = W^\prime + K$ to complete
the construction.

From the conditions above, it follows that $\{x_n\}$ is a Cauchy
sequence. Let $x = ||\cdot||$-$\lim x_n$. We have $x + K \su
\bigcap_{n=0}^{\infty} w$-${\rm cl}V_n$ and $x \in \p G_K$. Hence
$x+K\su {\rm cl}G$ and $(x + K)\cap F = \emptyset$. Since $K$ is a
compact set and $x \in \p G_K ,$ it follows that $(x + K) \cap \p G
\neq \emptyset$. The proof is complete. \finproof

Recall that a Banach space $X$ is called polyhedral \cite{K} if
the unit ball of any its finite-dimensional subspace is a
polytope. It was proved in \cite{F4} that a separable polyhedral
space admits a countable boundary. The next assertion shows that
fragmentable subsets of the unit sphere of a separable polyhedral
space are quite small.

\begin{corollary}\label{c1}
Let $G$ be an open bounded subset of a Banach space $X$, and let
$M\su\p G$ be a fragmentable set such that for any weak open set
$U$ with $U \cap M \neq \emptyset$, and for any weak closed set
$F$ with $F \cap G = \emptyset$, we have $(U \cap M)\setminus F
\neq \emptyset$. Then $G$ is not a weak $G_{\de}$ set. In
particular, if $X$ is polyhedral then, for every fragmentable set
$M \su S_X$, there is a weak open set $U$ which meets $M$ and a
finite number of hyperplanes $\{H_i\}_{i=1}^m$ in $X$, such that
$U \cap M \su \bigcup_{i=1}^m H_i .$
\end{corollary}

\n{\bf Proof.} We can assume that $0 \in G$. If we put $K = \{0\}$
and apply Proposition \ref{p1}, we see that $G$ is not a weak
$G_\delta$ set. If $X$ is polyhedral then \cite{F4} it has a
countable boundary and hence there is a sequence of hyperplanes
$\{H_i\}_{i=1}^\infty$ in $X$ with $S_X \su \bigcup_{i=1}^\infty
H_i$. Setting $F_n = \bigcup_{i=1}^n H_i$ for $n \geq 1$, using
the proof of Proposition \ref{p1}, we find a weak open set $U$
and $m \in \m{N}$ such that $(M\cap U)\setminus F_m =\emptyset$
and $M\cap U\not =\emptyset$. \finproof

\n{\bf Proof of Theorem \ref{t1}.} Let $M \su X$ be as in Definition
\ref{dk}. It will help to assume that $0 \in M$. If necessary, this
can be done by replacing $M$ with the set
$$
\left( \frac{M - z}{1 - ||z||} \right) \cap {\rm int}B_X
$$
where $z \in M$ is arbitrary. Assume that $G\su X$ is an open
bounded set and $0\in K\su G$, with $K$ a compact set which we
specify later. We will check the conditions of Proposition \ref{p1}.
First of all $0 \in M \cap G_K$. Now let $G_K , U, E,$ and $\ep
> 0,$ be as in Proposition \ref{p1}. Pick $x_0 \in U \cap M\cap G_K$
and by using the condition {\st}, find $\{x_i\}_{i=1}^n \su U \cap
M$ with $\|x_i -x_{i-1}\| < \ep$, $i = 1,\ldots,n$, $\|x_n\|\geq
1-\ep$. Assume that $x_n \in M\cap G_K$. Then since $||x_n|| \geq 1
- \ep$ and $x_n \in M\cap G_K \su G_K \su G \su B_X$, it follows
that $d(x_n,\p G) < \ep$. If $x_n \not\in M\cap G_K$ then there is
$m < n$ with $x_m \in M\cap G_K$ and $x_{m+1} \not\in M\cap G_K$.
Since $||x_m - x_{m+1}|| < \ep$, it follows that $d(x_m,\p G_K) <
\ep$. Set $y = x_m$ if $x_n \not\in M\cap G_K$, and $y = x_n$
otherwise. Hence $d(y, \p G_K)<\ep .$ Since $y \in M\cap G_K \su
G_K$ we get that $y\in G$. Having in mind that $E \cap G =
\emptyset$, we get $(y+K)\cap E=\emptyset .$

Now assume to contrary that $G$ is a weak $G_{\de}$ set. Put
$F=X\setminus G .$ Then $F$ is a weak $F_{\si}$ set and by
Proposition \ref{p1} there is $x\in X$ such that $x+K\su {\rm
cl}G,$ $(x+K)\cap F=\emptyset ,$ and $(x+K)\cap\p G\not
=\emptyset ,$ contradicting $\p G\su F.$

The proof of the second part of the theorem uses an idea from the
proof of \cite[Theorem 2]{F3}. Given a separable Banach space $X$
and a convex, bounded open set $G$ with $0\in G,$ we let $K =
T(B(\ell_2))$, where $T:\ell_2 \rightarrow X$ is a linear,
compact operator with dense range, and chosen so that $K$ is
contained in the interior of $G$. If $F$ is a weak $F_\sigma$ set
with $F \cap G = \emptyset$ then by Proposition \ref{p1} we obtain
$x \in \rm{ cl}G$ satisfying
\begin{equation}\label{eee}
x+K \su {\rm cl}G,\;\; (x+K) \cap \p G \neq\emptyset ,\;\;
(x+K)\cap \p G\cap F =\emptyset.
\end{equation}
Now assume to the contrary that $w^*$-exp $G^{\circ}$ is
countable. Then by Fact \ref{ff} the set ${\rm sm}({\rm cl}G)$ of
all smooth points of ${\rm cl}G$ is $w-F_{\si}.$ Put $F={\rm
sm}({\rm cl}G)$ and apply (\ref{eee}). We get a point $z\in
(x+K)\cap (\p G\setminus F).$ However by using that $K =
T(B(\ell_2))$ and ${\rm cl\: span}K=X$ it is easy to see that
$z\in F,$ a contradiction. The proof is complete. \finproof

\section{Examples}

Let $X$ be a Banach space with a normalized shrinking basis
$\{e_i\}$, such that there is a sequence of numbers $\{t_i\},\;
\lim_i t_i =0,$ with two further properties:
\begin{enumerate}
\item[(a)] $\sup_n ||\sum_{i=1}^n t_i e_i|| = \infty$;
\item[(b)] for any subsequence $\{t_{i_k}\}$ such that $\sup_n
||\sum_{k=1}^n t_{i_k}e_{i_k}|| < \infty$, the series
$\sum_{k=1}^{\infty} t_{i_k}e_{i_k}$ converges.
\end{enumerate}
We show there exists a relatively weakly compact subset $M \su B_X$,
satisfying condition {\st} of Theorem \ref{t1}.

Let $\{e_i^*\}$ be the biorthogonal sequence for $\{e_i\}$ and
$$
P_n x =\sum_{i=1}^n e_i^* (x)e_i ,\;\; x\in X,\;\; n=1,2,\ldots
$$
Denote
$$
M = \{x=\sum_{i\in\si}t_i e_i :\; \si\su \mathbb{N},\; |\si|<\infty,
\; ||P_n x||\leq 1,\; n=1,2,\ldots\}.
$$
Now pick $x_0 = \sum_{i \in \si_0}t_i e_i \in M$, $||x_0|| < 1 -
\ep$, and a weak open set $V$ containing $x_0$. Find $\de
> 0$ and $m\in \mathbb{N}$ such that
$$
x_0 \in U = \{u\in X:\; |e_i^* (x_0 -u)|<\de :\; i=1,\ldots,m\} \su
V.
$$
Given $\ep > 0$, find $l \in \mathbb{N}$ such that $|t_i| < \ep$ for
$i > l$. Denote $i_0 = \max\si$ and pick $j > \max\{i_0 ,l,m\}$. Set
$$
x_{k+1} = x_0 + \sum_{i=j}^{j+k}t_i e_i,\quad k = 0,1,\ldots
$$
Clearly, $\{x_k\} \su U$, $||x_k - x_{k+1}||< \ep$, $k =
0,1,\ldots$, and $\lim_k ||x_k|| = \infty$. Let $n$ be the minimal
index for which $||x_n|| < 1$ and $\|x_{n+1}\| \geq 1$. Then
$||x_k|| < 1$, $x_k \in M$, $k = 1,\ldots,n$, and $||x_n|| \geq
||x_{n+1}|| - ||x_n -x_{n+1}|| > 1 - \ep$.

Next we show that $M$ is relatively weakly compact. Given a sequence
$\{y_l\} \su M$, we have finite $\sigma_l \su \mathbb{N}$ such that
$y_l = \sum_{i \in \sigma_l} t_i e_i$ and $\sup_n ||P_n y_l||\leq 1$
for each $n$ and $l$. By taking a subsequence, we can find $\sigma
\su \mathbb{N}$ such that $\lim_l \sigma_l = \sigma$ in the
pointwise topology of the power set of $\mathbb{N}$. We enumerate
$\sigma$ as a strictly increasing sequence $\{i_k\}$. Clearly
$\lim_l P_{i_n} y_l = \sum_{k=1}^n t_{i_k} e_{i_k}$ for each $n$, so
by (b), $y = \sum_{k=1}^\infty t_{i_k} e_{i_k}$ converges in $X$.
Since $\{e_i\}$ is shrinking, it is evident that $w$-$\lim_l y_l =
y$.

\begin{example}\em{There is a separable Banach space $X$ with
shrinking basis which is $c_0$-saturated but does not contain a
bounded, open weak-$G_{\de}$ set. Moreover, for any equivalent norm
$|||\cdot|||$ on $X$, the set ${\rm exp}\,B_{(X,|||\cdot|||)^*}$ is
uncountable.}
\end{example}

Indeed, in \cite{L} a non-degenerate Orlicz function $M$ is
constructed such that there is a sequence $\{t_i\}$, $\lim_i t_i =
0$, with
$$
\sup_i \frac{M(Kt_i)}{M(t_i)}<\infty,
$$
for any $K > 0$, and
$$
\al_M = \sup\{q:\; \sup_{0 < \la, t \leq 1}\frac{M(\la
t)}{M(\la)t^q} < \infty\} = \infty.
$$
From \cite[p.\ 143]{LT}, it follows that the space $h_M$ is
$c_0$-saturated. By repeating some of the $t_i$ if necessary, we may
assume that $\sum_i M(t_i) = \infty$. Then the unit vector basis
$\{e_i\}$ of $h_M$ and the sequence $\{t_i\}$ satisfy the conditions
(a) and (b). Let us mention that in \cite{L}, it is shown that such
$h_M$ has no countable boundary for any equivalent norm.


\begin{thebibliography}{99}

\bibitem[B]{B} R. Bourgin, Geometric aspects of convex sets with the Radon-Nikod\'{y}m
property, Lecture Notes in Mathematics 993, Springer-Verlag, Berlin 1983.

\bibitem[F1]{F1} V. P. Fonf, On one property of Lindenstrauss-Phelps
spaces, {\em Funct. Anal. Appl.} 13 (1979) 66--67 (translation from
Russian).

\bibitem[F2]{F2} V. P. Fonf, Boundedly complete basic sequences,
$c_0$-subspaces, and injections of Banach spaces, {\em Israel. J.
Math.} 89 (1995) 173--188.

\bibitem[F3]{F3} V. P. Fonf, On exposed and smooth points of convex
bodies in Banach spaces, {\em Bull. London Math. Soc.} 28 (1996)
51--58.

\bibitem[F4]{F4} V. P. Fonf, Polyhedral Banach spaces, {\em Math. Notes
Acad. Sci. USSR} 30 (1981), 809--813 (translated from Russian).

\bibitem[FLP]{FLP} V.P. Fonf, J. Lindenstrauss, R.R. Phelps, Infinite
Dimensional Convexity, 599--670, {\em Handbook of the Geometry of
Banach Spaces}, Vol.1, Edited by W.B. Johnson and J. Lindenstrauss,
Elsevier Science, 2001.

\bibitem[K]{K} V. Klee, Polyhedral sections of convex bodies, {\em
Acta Math.} 103, (1960), 243--267.

\bibitem[L]{L} D. H. Leung, Some isomorphically polyhedral Orlicz
sequence spaces, {\em Israel J. Math.} 87 (1994), 117--128.

\bibitem[LP]{LP} J. Lindenstrauss, R. Phelps, Extreme point properties of convex bodies in
reflexive Banach spaces, {\em Israel. J. Math.} 6 (1968), 39--48.

\bibitem[LT]{LT} J. Lindenstrauss, L. Tzafriri, Classical Banach spaces I.
Springer-Verlag, Berlin-New York, 1977.

\bibitem[N]{N} I. Namioka, Radon-Nikod\'{y}m compact spaces and fragmentability,
{\em Mathematika} 34 (1987), 258--281.
\end{thebibliography}
\end{document}